\documentstyle[12pt,amssymb]{article}

\def\fnote#1#2{\begingroup\def\thefootnote{#1}\footnote{#2}\addtocounter{footnote}{-1}\endgroup}

\def\inbar{\vrule height1.5ex width.4pt depth0pt}
\def\IB{\relax{\rm I\kern-.18em B}}
\def\IC{\relax\,\hbox{$\inbar\kern-.3em{\rm C}$}}
\def\ID{\relax{\rm I\kern-.18em D}}
\def\IE{\relax{\rm I\kern-.18em E}}
\def\IF{\relax{\rm I\kern-.18em F}}
\def\IG{\relax\,\hbox{$\inbar\kern-.3em{\rm G}$}}
\def\IH{\relax{\rm I\kern-.18em H}}
\def\II{\relax{\rm I\kern-.18em I}}
\def\IK{\relax{\rm I\kern-.18em K}}
\def\IL{\relax{\rm I\kern-.18em L}}
\def\IM{\relax{\rm I\kern-.18em M}}
\def\IN{\relax{\rm I\kern-.18em N}}
\def\IO{\relax\,\hbox{$\inbar\kern-.3em{\rm O}$}}
\def\IP{\relax{\rm I\kern-.18em P}}
\def\IQ{\relax\,\hbox{$\inbar\kern-.3em{\rm Q}$}}
\def\IR{\relax{\rm I\kern-.18em R}}
\def\IT{\relax{\rm I\kern-.18em T}}
\def\ZZ{\relax{\sf Z\kern-.4em Z}}

\def\fnote#1#2{\begingroup\def\thefootnote{#1}\footnote{#2}\addtocounter
{footnote}{-1}\endgroup}
\def\beq{\begin{equation}}
\def\eeq{\end{equation}}
\def\bea{\begin{eqnarray}}
\def\eea{\end{eqnarray}}

\def\notin{\ \hbox{{$\in$}\kern-.51em\hbox{/}}}

  \def\E1Fq{E_1/\IF_q}

\def\notdiv{{\relax{~|\kern-.34em /~}}}

\begin{document}
\parindent=0pt
\hfill NSF$-$ITP$-$02$-$99

\vskip 0.7truein

\centerline{\large {\bf ASPECTS OF CONFORMAL FIELD THEORY}}

\vskip .1truein

\centerline{\large {\bf FROM CALABI-YAU
ARITHMETIC}\fnote{\diamond}{To appear in the Proceedings of the
{\it Workshop on Arithmetic, Geometry and Physics around
Calabi-Yau Varieties and Mirror Symmetry}, eds. J. Lewis and N.
Yui, American Math. Society.}}

\vskip 0.3truein

\centerline{\sc Rolf Schimmrigk\fnote{\dagger}{email:
netahu@yahoo.com}\fnote{+}{This work was supported in part by NSF
Grant \# PHY99-07949. The author is grateful to the Kavli
Institute for Theoretical Physics at Santa Barbara for hospitality
as a Scholar.}}

\vskip .2truein

\centerline{\it Kennesaw State University} \vskip .05truein
\centerline{\it 1000 Chastain Rd, Kennesaw, GA 30144}

\vskip .6truein

\baselineskip=19pt

\centerline{\bf ABSTRACT:} \vskip .2truein
This paper describes a
framework in which techniques from
 arithmetic algebraic geometry are used to formulate a direct
 and intrinsic link between the geometry of Calabi-Yau
 manifolds
 and aspects of the underlying conformal field theory.
 As an application the algebraic number field
 determined by the fusion rules of the conformal field theory
 is derived from the number theoretic structure of
 the cohomological Hasse-Weil L-function determined by
  Artin's congruent zeta function of the algebraic variety.
In this context a natural number theoretic characterization arises
for the quantum dimensions in this geometrically determined
algebraic number field.

\vfill

{\sc PACS Numbers and Keywords:} \hfill \break Math:  11G25
Varieties over finite fields; 11G40 L-functions; 14G10  Zeta
functions; 14G40 Arithmetic Varieties \hfill \break Phys: 11.25.-w
Fundamental strings; 11.25.Hf Conformal Field Theory; 11.25.Mj
Compactification

\renewcommand\thepage{}
\newpage

\baselineskip=22pt
\parskip=.1truein
\parindent=0pt
\pagenumbering{arabic}

\section{Introduction}

One of the results to emerge from the second period of string
theory  has been the discovery of mirror symmetry some ten years
ago. The first large scale cohomological evidence of mirror
symmetry presented in \cite{cls90} was soon understood to be a
consequence of a mirror construction that mapped Calabi-Yau
orbifolds into complete intersection varieties \cite{ls90b}.
Contemporaneously, a much more precise way of constructing mirror
pairs by orbifolding exactly solvable string compactifications was
discovered in ref. \cite{gp90}. It is the conformal field theory
picture which gives confidence to the conclusion that
topologically rather different mirror pairs of Calabi-Yau
varieties should define isomorphic compactifications of string
theory. This insight has been used by physicists and
mathematicians to derive a number of interesting consequences,
such as the computation of the number of rational curves on
Calabi-Yau varieties \cite{cdgp90}. A review of many of the
results can be found in \cite{ck}.

At present, mirror symmetry, in all its geometric ramifications,
is still poorly understood, and the conformal field theoretic
constructions remain an important tool in its exploration. The
simplest class of models in which conformal field theory has been
linked to Calabi-Yau varieties is provided by Gepner's
construction of tensor models of minimal $N=2$ superconformal
field theories \cite{g87}. When Gepner originally constructed
non-geometric string compactifications with 2-dimensional
interacting exactly solvable theories on the worldsheet his
expectation was to have obtained an entirely new class of
consistent string backgrounds. But when comparing the spectra of
several of his models with the Hodge numbers of some existing
Calabi-Yau manifolds \cite{chsw85,s87} he noticed that,
surprisingly, there was agreement between the spectra of these
rather different looking models \cite{g87b}. Further work
established a more systematic relation, initially based on a
Landau-Ginzburg mean field theory description of the conformal
field theory \cite{m89,vw89,lvw89,ls90a,fkss90}, and later on a
sigma-theoretic analysis \cite{w92}.

In the Landau-Ginzburg picture part of the structure of the
conformal field theory is encoded in the chiral ring, determined
by the superpotential, which for diagonal minimal models takes the
form $W(\Phi_i) = \sum_{i=0}^s \Phi_i^{n_i}(z,{\bar z})$, where
the $\Phi_i(z,{\bar z})$ are chiral primary fields on the string
world sheet Riemann surface. In the sigma-theoretic picture both
the Landau-Ginzburg theory and the Calabi-Yau variety are
recovered as different phases of the model.

Both the Landau-Ginzburg and the sigma-theoretic approach
illuminate the relation between Calabi-Yaus and conformal field
theories. What is lacking, however, is a direct, intrinsic, and
mathematically rigorous framework which allows to establish a link
between algebraic varieties and conformal field theories. One
might hope for a framework in which it is possible to derive the
essential ingredients of the conformal field theory directly from
the algebraic variety itself, without the intermediate
Landau-Ginzburg formulation. A priori it might appear unlikely
that such an approach exists because numbers associated to
algebraic varieties usually are integers (such as dimension,
cohomology dimensions and indices associated to complexes),
whereas the numbers appearing in the underlying exactly solvable
conformal field theory usually are rational numbers (such as the
central charge and the anomalous dimensions).

This article describes a strategy to formulate an intrinsic,
direct, and mathematically rigorous framework which allows to
derive certain conformal field theoretic quantities directly from
the algebraic Calabi-Yau variety \cite{rs01}. The idea is to use
concepts from
 arithmetic algebraic geometry, in which algebraic varieties $X$
 are defined not over a continuous field, like the real numbers
 ${\mathbb R}$,
 or the complex number field ${\mathbb C}$, but over discrete
 finite fields
 ${\mathbb F}_q$, where $q \in {\mathbb N}$ denotes the number of
 elements of the field.
The particular field of reduction of the variety is specified by
writing $X/{\mathbb F}_q$, leading to a reduced variety consisting
of a finite number of points.

In the past it has not been particularly useful in physics to
consider manifolds over fields other than the complex or real
numbers. All work done in string compactification over the last
fifteen years has been done over the complex numbers\footnote{A
recent exception appears in these proceedings
\cite{cov00,cov01}.}. Lattice constructions usually are viewed as
a preliminary step to taking some kind of continuum limit which is
considered to describe the correct model. From a physical point of
view the choice of any finite $q$ in the reduction of a variety
$X/{\mathbb F}_q$ would appear to be arbitrary and ill-motivated.
For small $q$ the field ${\mathbb F}_q$ would define a large scale
lattice structure which one might expect to provide only rough
information about the structure of the variety. More sensible is
to consider an infinite sequence of ever larger finite fields
which probe the variety at ever smaller scales. It is this
consideration which leads us to the concept of counting the number
of solutions $N_{r,p}(X)=\#(X/{\mathbb F}_{p^r})$ and to ask what,
if any, interesting information is provided by the numbers
determined by the extensions ${\mathbb F}_{p^r}$ of degree $r$ of
the finite field ${\mathbb F}_q$. The following paragraphs outline
the strategy envisioned in this program and briefly introduce the
arithmetic ingredients used in this paper.

The starting point of the arithmetic considerations in Section 4
is to arrange the sequence of reductions of the variety over the
finite fields ${\mathbb F}_{p^r}$ into a useful form. This can be
achieved via Artin's congruent zeta function, essentially defining
an exponential sum of a generating function constructed from the
numbers $N_{r,p}(X)$. For reasons just described we should not
restrict the construction to local considerations at particular
integers. Hence we need some way to pass to a global description
in which all reference to fixed scales has been erased. This leads
in Section 5 to the concept of the global Hasse-Weil L-function of
a Calabi-Yau variety. This L-function will collect the information
at all rational primes. This step can be taken because the Weil
conjectures proven by Dwork and Deligne show that in the present
context Artin's zeta function is a rational function determined by
the cohomology of the variety. The Hasse-Weil L-function therefore
is a cohomological L-function.

At this point the only number theoretic algebraic structures that
have appeared are the finite fields ${\mathbb F}_p$ and their
extensions. A brief review of some of the aspects of the
underlying conformal field theoretic aspects, provided in Section
2, reveals that this is not enough. As mentioned above, the
structures that enter in Gepner's construction and its
generalizations are rational, i.e. they lead to central charges
and spectra of anomalous dimensions which are rational numbers. It
is not clear how to recover these numbers from the intrinsic
geometry of the variety. It turns out that it is more useful in
the present context to encode the conformal field theoretic
information in an alternative way by mapping the central charge
and the anomalous dimensions via the Rogers dilogarithm into the
quantum dimensions associated to the physical fields. Section 3
reviews that these (generalized) quantum dimensions are elements
of certain algebraic number fields which are determined by the
fusion rules of the conformal field theory.

This suggests that we consider the structure of algebraic number
fields in more detail. It is in this context that Hecke introduced
a general notion of L-functions which are determined by the prime
ideals of the number field, generalizing earlier results obtained
by Dirichlet and Dedekind. Such number field L-functions, reviewed
in Section 7, essentially count the inverse of the norms of ideals
in the number field weighted by a character of this field.

The surprise, for the uninitiated, is that the cohomological
Hasse-Weil type L-function of the Calabi-Yau variety in fact
happens to be determined by number field L-functions of the type
introduced by Hecke. This is the link that allows us to recover in
Sections 7 and 8 the number theoretic framework relevant to the
underlying conformal field theory directly from the intrinsic
geometry of the variety. Once the fusion field has been identified
one can then further explore the arithmetic r\^{o}le played by the
quantum dimensions in this field.

The strategy used here can be summarized as follows. First
consider the arithmetic structure of Calabi-Yau varieties via
Artin's congruent zeta function. Next define the global
cohomological Hasse-Weil L-function via local factors from the
congruent zeta function. Finally interpret the Hasse-Weil L-series
as Hecke L-series of an algebraic number field. In this way we can
recover the algebraic fusion field of the conformal field theory
from the algebraic variety.

The program to use arithmetic geometry to illuminate the conformal
field theoretic structure of Calabi-Yau varieties originated
several years ago. The idea initially was to develop further some
results derived by Bloch and Schoen in support of the
Beilinson-Bloch conjectures, and apply them in the context of the
conformal field theory/Calabi-Yau relation \cite{s95}. The final
section of the paper contains a brief outline of the concepts
involved in framework. More recently arithmetic considerations
have been discussed in different contexts in refs. \cite{m98} and
\cite{cov00}.

\section{Conformal Field Theory}

\subsection{Rational Theories}

The simplest constructions of exactly solvable Calabi-Yau
varieties involve conformal field theories with special features,
described by rational exactly solvable field theories with N=2
supersymmetry.

These are defined by infinite dimensional Lie-algebras defined by
the so-called Virasoro algebra defined in terms of the Fourier
modes of the energy momentum tensor on the 2-dimensional theory
$$T(z) = \sum_{n \in {\mathbb Z}} \frac{L_n}{z^{n+2}}$$ for which
the algebra takes the form $$ [L_m, L_n] = (m-n) L_{m+n} +
\frac{c}{12} m(m^2-1) \delta_{m+n,0}. $$

N=2 supersymmetry requires that the Virasoro algebra is extended
by some affine Lie algebra currents $$J(z) = \sum_{n\in {\mathbb
Z}} \frac{J_n}{z^{n+1}},$$
 e.g. ${\rm SU}(n)_k$ at level $k$, described by the commutator
 relation of its components $J_n = J_n^a T_a$
$$
[J^a_m, J^b_n] = i f^{ab}_c J^c_{m+n} + km
\delta_{m+n,0}\delta^{ab}.
$$

Rationality of the theory means that the Hilbert space decomposes
into a finite number of blocks
$$
{\mathcal H} = \bigoplus_{\ell=1}^N {\mathcal H}_{\ell}.
$$

Such models are characterized by a set of rational numbers: \hfill
\break $\bullet$ The number $k$ appearing in the defining algebra
is the level of the conformal field theory and takes values in the
positive rational integers $k \in {\mathbb N}$. \hfill \break
$\bullet$ The central charge $c$ in the Virasoro algebra takes
rational values $c\in {\mathbb Q}$. \break $\bullet$ The anomalous
dimensions $\Delta_{\ell}$ of the physical (chiral primary) fields
$\phi_{\ell}(z,{\bar z})$ on the worldsheet surfaces swept out by
the string (torus) are determined to take rational values as well
$\Delta_{\ell}\in {\mathbb Q}$.

The classification of all rational conformal field theories with
N=2 worldsheet supersymmetry has not been completed yet. The most
notable constructions that have been analyzed so far are based on
affine coset constructions, in particular the generalizations of
Gepner models introduced by Kazama and Suzuki \cite{ks89}. In the
following the focus will be on the simpler Gepner models.

\subsection{SU(2)$_k$ Building Block}

The simplest examples of Calabi-Yau manifolds that lead to exactly
solvable models are described by extensions of the Virasoro
algebra defined by currents of the SU(2) affine theory at some
level $k$, for which the central charge is given by
$$c=\frac{3k}{k+2},$$ and the anomalous dimensions are determined
to be
$$ \Delta_{(k)}^{\ell} = \frac{\ell (\ell+2)}{4(k+2)},~~~~~~~
\ell=0,\dots, k.$$

Associated to the anomalous dimensions are certain conformal field
theory characters which are functions on the worldsheet of the
string (genus 1 curve). More precisely, consider maps on the upper
half-plane ${\mathfrak H}$ $$ \chi^{\ell}: {\mathfrak H} \times
{\mathbb C}^2 \longrightarrow {\mathbb C}$$ defined as $$
\chi^{\ell}(\tau, z, u):= e^{-2\pi i ku} {\rm tr}_{{\mathcal
H}_{\ell}} q^{L_0 - \frac{c}{24}} e^{2\pi i uJ_0}.$$

These characters can be expressed via the string functions
$$c^{\ell}_m(\tau) = \frac{1}{\eta^3(\tau)}
\sum_{\stackrel{\stackrel{-|x|<y\leq |x|}{(x,y)~{\rm
or}~(\frac{1}{2}-x,\frac{1}{2}+y)}}{\in {\mathbb Z}^2 +
\left(\frac{\ell+1}{2(k+2)},\frac{m}{2k}\right)}} {\rm sign}(x)
e^{2\pi i \tau((k+2)x^2-ky^2)}$$ and the classical theta functions
$\Theta_{m,k}$ $$ \Theta_{n,m}(\tau,z,u) = e^{-2\pi i m u}
\sum_{\ell \in {\mathbb Z} + \frac{n}{2m}} e^{2\pi i m \ell^2 \tau
+ 2\pi i \ell z}$$ as
\begin{equation} \chi^{\ell}(\tau, z,u) =
\sum_{\stackrel{m=-k+1}{\ell=m~{\rm mod}~2}}^k c^{\ell}_m (\tau)
\Theta_{m,k}(\tau, z,u).\end{equation}

This set of characters is closed under modular transformations.
For the translation generator $\tau \mapsto \tau +1$ the
characters acquire a simple scaling factor
$$ \chi^{\ell} \left(\tau+1,z,u\right) =
e^{2\pi i\Delta_{\ell}} \chi^{\ell}(\tau,z, u)$$ while the
S-generator $\tau \mapsto -1/\tau$ leads to the relation
$$ \chi^{\ell} \left(-\frac{1}{\tau}, \frac{z}{\tau},
u+\frac{z^2}{2\tau}\right) = e^{\pi i kz^2/2} \sum_{m=0}^k S_{\ell
m} \chi^m(\tau, u),$$ with modular $S$-matrix
\begin{equation}
S_{\ell m} =
\sqrt{\frac{2}{k+2}}~~\sin\left(\frac{(\ell+1)(m+1)\pi}{k+2}\right)
\label{S-matrix} \end{equation}

\subsection{N=2 Supersymmetric models}

The simplest class of N=2 supersymmetric exactly solvable theories
is built in terms of the affine SU(2) theory as a coset model
\begin{equation}  \frac{{\rm SU(2)}_k \otimes {\rm U(1)}_2}{{\rm
U(1)}_{k+2,{\rm diag}}}.\end{equation} Coset theories $G/H$ lead
to central charges of the form $c_G - c_H$, hence the
supersymmetric affine theory at level $k$  still has central
charge $c_k=3k/(k+2)$. The spectrum of anomalous dimensions
$\Delta^{\ell}_{q,s}$ and U(1)$-$charges $Q^{\ell}$ of the primary
fields $\Phi^{\ell}_{q,s}$ at level $k$ is given by
\begin{eqnarray} \Delta^{\ell}_{q,s} &=&
\frac{\ell (\ell +2)-q^2}{4(k+2)} + \frac{s^2}{8} \nonumber \\
Q^{\ell}_{q,s} &=& \frac{q}{k+2} - \frac{s}{2}, \end{eqnarray}
where $\ell\in \{0,1,\dots,k\}$, $\ell+q+s \in 2{\mathbb Z}$, and
$|q-s|\leq \ell$. Associated to the primary fields are characters
defined as
\begin{equation}
\chi^{\ell}_{q,s}(\tau, z,u) = e^{-2\pi i u} {\rm tr}_{{\mathcal
H}^{\ell}_{q,s}} e^{2\pi i\tau (L_0 -\frac{c}{24})} e^{2\pi i
J_0}, \end{equation} where the trace is over the projection of
${\mathcal H}^{\ell}_{q,s}$ to definite fermion number mod 2 of a
highest weight representation of the N$=$2 superfonformal algebra.
It is of advantage to express these maps in terms of the string
functions and theta functions, leading to the form
\begin{equation} \chi^{\ell}_{q,s} = \sum
c^{\ell}_{q+4j-s}\Theta_{2q+(4j-s)(k+2),2k(k+2)}, \end{equation}
because it follows from this representation that the modular
behavior of the $N=2$ characters decomposes into a product of the
affine SU(2) structure in the $\ell$ index and into
$\Theta$-function behavior in the charge and sector index. More
precisely, the modular behavior of these characters is given by
\begin{equation} \chi^{\ell}_{q,s}(\tau +1,z,u) = e^{2\pi i
\Delta^{\ell}_{q,s}} \chi^{\ell}_{q,s}(\tau, z,u)
\end{equation} for the modular shift $\tau \mapsto \tau +1$, and
by
\begin{eqnarray}
&&\chi^{\ell}_{q,s}\left(-\frac{1}{\tau}, \frac{z}{\tau},
u+\frac{z^2}{2\tau}\right) =  \nonumber \\
& & \phantom{w} \frac{1}{\sqrt{2}(k+2)}
\sum_{\stackrel{\ell',q',s'}{\ell'+q'+s'=0~{\rm mod}~2}} e^{\pi i
\frac{qq'}{k+2}} e^{-\pi i \frac{ss'}{2}} \sin\left(\pi
\frac{(\ell+1)(\ell'+1)}{k+2}\right)
\chi^{\ell'}_{q',s'}(\tau,z,u) \nonumber \\ \end{eqnarray} for the
generator $\tau \mapsto -1/\tau$.

Because the transformation behavior of these characters splits
into the nontrivial interacting SU(2) theory and a pair of U(1)
theories, the partition function of these minimal $N=2$ models can
be written as \begin{equation} Z_{N=2}(\tau, {\bar \tau}) =
\frac{1}{2} \sum_{\stackrel{\ell, {\bar \ell}}{q,{\bar q},s, {\bar
s}}} A_{\ell, {\bar \ell}} M^{q,{\bar q}} M^{s, {\bar s}}
\chi^{\ell}_{q,s}(\tau) \chi^{\bar \ell}_{{\bar q},{\bar s}}({\bar
\tau}).\end{equation}

\subsection{String Compactification with Gepner Models}

String theory in a supersymmetric geometric background lives in
ten dimensions because of anomaly cancellation constraints.
Decompactification of four of the initially ten compact dimensions
leaves a compact six-dimensional variety. In the context of
theories with $N=1$ spacetime supersymmetry each of these compact
dimensions contributes an anomaly of 3/2. This means that in order
to construct an internal theory based on affine SU(2) theories one
has to glue several of these models together in a tensor product
structure to provide the correct conformal anomaly. The number of
SU(2)$_{k_i}$-factors one needs to tensor to make
 a Calabi-Yau variety is determined by the sum
$$ c_{\rm tot} = \sum_{i=1}^r c_i =
\sum_{i=1}^r \frac{3k_i}{k_i+2}=9.$$

 The physical spectrum of the massless modes of the compactified
theory is in part described by the cohomology of the Calabi-Yau
variety. In the conformal field these fields are determined by the
two-dimensional chiral primary fields and the dimension of the
geometric cohomology theories can be computed by combinatorics of
the anomalous dimension of the conformal field theory. More
precisely, the multiplicities are determined by the partition
function of the theory, which takes the form (in the light cone
gauge)

\begin{equation}
Z_{het} = \mid ({\rm Im}~ \tau)^{-1}\eta(\tau)^{-4}\mid
\sum_{\lambda \in \{o,v,s,\bar {s}\}}
   B_{\lambda}^{SO(2)}(\tau)
   B_{\bar {\lambda}}^{E_8\times SO(10)}(\bar {\tau})
   Z_{\rm int}(\tau,\bar {\tau}),
\end{equation}

\noindent where the 4$d$--part comes from the spacetime
variables.
The $SO(2)$ characters $B_{\lambda}^{SO(2)}$ come from the
spacetime fermions, whereas the characters $B_{\bar
{\lambda}}^{E_8\times SO(10)}$ are determined by the internal
fermions, and  the internal partition function $Z_{\rm int}$ is
determined by the tensor models and will be given by products of
the form $\prod_{i=1}^r \chi^{\ell_i}_{q_i,s_i}$. Collecting the
affine indices into an $\vec \ell$--index and the charge and the
sector indices $(q_i,s_i)$ as well as the external index (coming
from the level 2 theta function of the 4--d fermions) into one
$(2r+1)$ vector, we can write the partition function (modulo the
external transverse bosons) as
\begin{equation}
\chi^{\vec \ell}_{\vec \mu} = B_{\mu_0}^G
\prod_{i=1}^r\chi^{\ell_i}_{\mu_i,\mu_{r+i}},
\end{equation}
\noindent where $G$ stands for SO(2) in the supersymmetric sector,
and for SO(26) or $E_8 \times SO(10)$ in the Kac--Moody sector.
The index $\mu_0$ labels the highest weight of a given
representation of $G$. The standard fermionic invariant for the
theta--function products then leads to
\begin{equation}
Z \sim \sum_{\vec \ell,\vec {\bar \ell}} A_{\vec \ell,\vec {\bar
\ell}}
       \sum_{\mu,m_i,n_i} (-1)^{m_0 + n_0}
       \chi^{\vec \ell}_{\vec {\mu} + m_i \vec {\beta^i}}
       \chi^{\vec {\bar \ell} *}_{\vec {\mu'} + n_i \vec {\beta^i}}
\label{int-part}\end{equation} for the compactified theory.
Consider eq. (\ref{int-part}) in more detail. The multi-matrix
$A_{\vec \ell,\vec {\bar \ell}} = \prod_{i=1}^r A_{\ell_i,{\bar
\ell}_i}$ represents a choice of the affine invariant for each of
the individual factors in the internal sector. The
$(2r+1)$--vectors $\vec \beta^i, ~\vec \mu,~\vec \mu'$ are chosen
in such a way as to implement requirements (B) and (C). To that
account define $\vec {\beta^0} = (\bar s|1,\dots,1)$ and $\vec
{\beta^i} = (v|0,\dots,0,2,0,\dots,0),~~i=1,\dots,r$ where 2 is in
the $(r+1+i)^{th}$ entry and $\bar s, v$ are highest weights for
the group $G$. The vectors $\vec \mu$ are then constrained by
\begin{equation}
\beta^i \cdot \mu \in 2{\mathbb Z} + \delta^{i0}
\end{equation}
\noindent where the scalar product is defined as follows. Recall
that the charge indices $q_i$ belong to level $(k+2)$ theta
functions whereas the external and the sector indices belong to
level 2 theta functions. It is then natural to define
\begin{equation}
\vec \beta^i \cdot \mu = \frac{\beta_0^i \mu_0}{2} +
 \sum_{j=1}^r \left(-\frac{\beta_j^i \mu_j}{k+2}
    +\frac{\beta_{r+j}^i \mu_{r+j}}{2}\right).
\end{equation}
\noindent For $i=0$ this condition implements the requirement
that
the U(1) charge must be quantized, while for $i\neq 0$ the
requirement of appropriate sector coupling (Neveu-Schwarz to
Neveu-Schwarz and Ramond to Ramond) is implemented in the
partition function. The integers $m_i$ and $n_i$ run over the
integers modulo $\alpha$ where $\alpha$ is such that $\beta_i
\alpha /2(k_i +2)$ are integers for all $i$. The vectors $~\vec
\mu,~\vec \mu'$ only differ in their external entry in such a
way
as to implement the transition from the superstring to the
heterotic string. This is accomplished by permuting the highest
weights $(o,v,s,{\bar s})$ of the gauge group $G$ into
$(v,o,s,{\bar s})$.

 The physical couplings in the conformal field theory are
determined mostly (not completely) by the anomalous dimensions of
the fields involved.

\subsection{An example: the quintic threefold}

An example of an exactly solvable Calabi-Yau threefold is the
quintic, described as the Fermat type variety in ordinary
projective 4-space ${\mathbb P}_4$
\begin{equation} X_5=  \left\{p(z_i) = z_0^5 + \cdots
+ z_4^5 = 0 \right\} \subset {\mathbb P}_4.\end{equation}

In terms of the intermediate Landau-Ginzburg formulation the
polynomial $p(z_i)$ originates from the superpotential
\begin{equation}
W(\Phi_i(x,{\bar x})) = \Phi_0^5(x,{\bar x}) + \cdots +
\Phi_4^5(x,{\bar x})
\end{equation}
in terms of the chiral primary super fields
$\Phi_i(x,{\bar x}) =
z_i + \psi_i \theta + \cdots$ on the worldsheet parametrized by
variables $(x,{\bar x})$. The complex coordinates $z_i$ of the
projective variety are obtained by taking the size of the string
to zero.

Conformal field theoretically we need theories of level $k=d-2$,
where in general $d$ is the degree of the superpotential. In the
present case this leads to SU(2) models at level 3, hence a
tensor
product of five theories is needed $X_5 \cong {\rm
SU(2)}_3^{\otimes 5}$. The cohomology of the projective variety
can be reconstructed in terms of primary fields such that the
total anomalous dimension $\Delta_{\rm tot} = \sum_{\ell}
\Delta^{\ell}_{\ell,0} =1$ and the U(1) charges are $\pm 1$.

\section{Number Fields from Rational Conformal Field Theories.}

The idea in the following is to map part of the conformal field
theory data into an interesting number field which can be
recovered in a geometric way from the variety.

The starting point is the collection of conformal field theory
characters associated to the anomalous dimensions of the primary
fields of the theory, more precisely their transformation behavior
with respect to the modular transformations on the worldsheet.
These transformations involve in particular the generator $S: \tau
\mapsto -1/\tau$. The nontrivial behavior of the N=2 characters is
carried by their affine SU(2) structure and involves the modular
$S$-matrix (\ref{S-matrix}). (The $T$-matrix associated to $\tau
\mapsto \tau +1)$ does not enter in the present context).

With these S-matrices one can define so-called quantum dimensions
$$Q_{\ell} = \frac{S_{0\ell}}{S_{00}}$$ and their generalized
cousins
$Q_{\ell m} = \frac{S_{\ell m}}{S_{0m}}$ for the SU(2) theory at
level $k$.

These numbers are of interest because it is possible to
reconstruct the anomalous dimensions from them.

Define the Rogers dilogarithm function as $$L(z) = Li_2(z) +
\frac{1}{2} \log(z) \log(1-z),$$ where $$Li_2(z) = \sum_{n\in
{\mathbb N}} \frac{z^n}{n^2}$$ is Euler's dilogarithm. This
function leads to transcendental values when evaluated at
algebraic numbers. It is therefore surprising that it can be used
to express rational numbers. This however is precisely what has
emerged.

{\bf Theorem.}[Kirillov-Reshetikhin, Lewin]{\it Let $Q_i =
S_{0i}/S_{00}$ be the quantum dimensions of the affine Lie algebra
SU(2)$_k$ at level $k$. Then}
\begin{equation} \frac{1}{L(1)} \sum_{\ell=1}^k
L\left(\frac{1}{Q_{\ell}^2}\right) = \frac{3k}{k+2}.
\end{equation}

Even more surprising is the result of Kuniba and Nakanishi
according to which the anomalous dimensions $\Delta_j^{(k)}$  of
the primary fields of the affine SU(2) theory at level $k$ can be
reconstructed from the quantum dimensions
\begin{equation}
\frac{1}{L(1)} \sum_{\ell=1}^k L\left(\frac{1}{Q_{\ell
m}^2}\right) = \frac{3k}{k+2} - 24 \Delta_{(k)}^m +6m.
\end{equation}

These relations provide an alternative way of parametrizing the
anomalous spectrum of fields at fixed central charge via the
quantum dimensions. The reason this is of interest in the present
context is because the quantum dimensions lead to an interesting
field.

{\bf Theorem.}[Boer \& Goeree \cite{dg91}] {\it The field
extension $K/{\mathbb Q}$ defined by the quantum dimensions of the
affine theory SU(2)$_k$ is normal with abelian Galois group ${\rm
Gal}(K/{\mathbb Q})$.}

Classical results by Kronecker and Weber then immediately
determine the type of field that contains the quantum dimensions.

{\bf Theorem.}[Kronecker-Weber]{\it The quantum dimension field
extension is contained in a cyclotomic extension ${\mathbb
Q}(\mu_m)/{\mathbb Q}$ for some integer $m$.}

This follows from the fusion rules for the chiral primary fields
$$ \phi_{\ell} \star \phi_m = \sum_m N_{\ell m}^r \phi_r,$$
leading to an additive behavior of the anomalous dimensions
$\Delta_{\ell}$ of $\phi_{\ell}$, and E. Verlinde's result for the
fusion matrices \cite{ev88}
$$N_{\ell m}^n = \sum_m \frac{S_{\ell r}S_{mr}}{S_{0r}}
(S^{-1})_{nr}.$$

The alternative parametrization of the anomalous dimensions in
terms of the quantum dimensions makes it possible to rephrase the
problem of an intrinsic derivation of these characteristics of the
conformal field theory as a derivation of the field of quantum
dimensions directly from the intrinsic geometry of the variety.

\section{The Congruent Zeta Function}

The basic idea is to derive the arithmetic structure of the
underlying conformal field theory directly from the arithmetic
structure of the variety $X$. The complete diophantine information
of $X$ is contained in the cardinalities of the reduced variety
$N_{r,p}(X) = \#(X/{\mathbb F}_{p^r})$. These numbers can be used
to define a generating function which, in analogy to Dedekind's
zeta function of an arbitrary algebraic number field, was
introduced in the context for curves over finite fields by E.
Artin in 1924 \cite{a24} as
\begin{equation}
Z(X/{\mathbb F}_p, t) = \exp\left(\sum_{r\in {\mathbb N}} \#
\left(X/{\mathbb F}_{p^r}\right) \frac{t^r}{r}\right).
\label{congzeta} \end{equation}

The motivation for arranging the number of solutions of
$X/{\mathbb F}_{p^r}$ into the form (\ref{congzeta}) is that this
expression leads to rational functions in the formal variable $t$.
This was first shown for curves $C_g$ of arbitrary genus $g$ by
F.K. Schmidt via the Riemann-Roch theorem \cite{fks31}\cite{hh34}
\begin{equation}
Z(C_g/{\mathbb F}_p,t) = \frac{P^{(p)}(t)}{(1-t)(1-pt)} ,
\end{equation} where $P^{(p)}(t)$ is a polynomial of degree
$2g$.

After much work by Artin, Hasse, Deuring and Weil, culminating in
the $p$-adic analysis of Dwork \cite{d60} and Deligne \cite{d74},
it was shown that more generally the zeta function of
$d$-dimensional smooth projective varieties over finite fields has
a number of important structural properties which justify its
definition.

{\bf 1.} The zeta function of a variety $X$ of complex dimension
$d$ is a rational function in the variable $t$
\begin{equation}
Z(X/{\mathbb F}_p,t)=\frac{\prod_{j=1}^d
P^{(p)}_{2j-1}(t)}{\prod_{j=0}^d P^{(p)}_{2j}(t)},\end{equation}
where the $P_i^{(p)}(t)$ are polynomials whose degrees are
determined by the Betti numbers $b_i={\rm dim~ H}^i_{\rm
DeRham}(X)$:
\begin{eqnarray} P_0^{(p)}(t)&=&1-t \nonumber \\
        P_{2d}^{(p)}(t)&=& 1-p^dt \nonumber
\end{eqnarray}
and for $1\leq i \leq 2d-1$
$$P_i^{(p)}(t) = \prod_{j=1}^{b_i}
\left(1-\beta^{(i)}_j(p) t\right),$$
 with algebraic integers $\beta^{(i)}_j(p)$.

{\bf 2.} $Z(X/{\mathbb F}_p,t)$ satisfies a functional equation
\cite{k94}
\begin{equation} Z\left(X/{{\mathbb F}_p}, \frac{1}{p^dt} \right)
= (-1)^{\chi +\mu} p^{d\chi/2} t^{\chi} Z(X/{{\mathbb F}_p} , t),
\end{equation} where $\chi$ is the Euler number of the variety $X$
over the complex numbers ${\mathbb C}$. Furthermore, $\mu$ is zero
when the dimension $d$ of the variety is odd, and $\mu$ is the
multiplicity of $-p^{d/2}$ as an eigenvalue of the action induced
on the cohomology by the Frobenius automorphism $\Phi:
X\longrightarrow  X,~x\mapsto x^p$.

{\bf 3.} The norms of the algebraic integers $\beta^{(i)}_j(p)$
satisfy the Riemann hypothesis \cite{d74} \begin{equation}
|\beta^{(i)}_j(p)| = p^{i/2},~~~~\forall i.\end{equation}

{\bf 4.1} There are four cases of interest:

({\bf a.}) Calabi-Yau curves (Elliptic curves) $$Z(E/{\mathbb
F}_p, t) = \frac{P_1^{(p)}(t)}{(1-t)(1-pt)}$$ with quadratic
$P_1^{(p)}(t)$. Such curves are used for compactifications of
various string theories to eight dimensions and provide a simple
framework in which different types of duality conjectures can be
investigated in great detail.

({\bf b.}) Calabi-Yau two-folds with finite fundamental group (K3
surfaces) have a Hodge diamond of the form $$\matrix{&   &1  &  &
\cr
          &0  &   &0 & \cr
       1  &   &20 &  &1, \cr }$$
leading to a congruent zeta function $$Z({\rm K3}/{\mathbb F}_p,t)
= \frac{1}{(1-t)P_2^{(p)}(t)(1-p^2t)}$$ with $P_2^{(p)}$ of degree
22.

({\bf c.}) Calabi--Yau threefolds with finite fundamental group
lead to Hodge diamonds of the form $${\small \matrix{
                     &   &         &1       &         &  & \cr
                     &   &0        &        &0        &   &\cr
                     &0  &         &h^{1,1} &         &   &\cr
                     1&   &h^{2,1}  &        &h^{2,1}  &
                     &1,\cr}}
             $$
resulting in zeta functions of the form $$ Z(X/{\mathbb F}_p, t)=
\frac{P^{(p)}_3(t)}{(1-t)P^{(p)}_2(t)P^{(p)}_4(t)(1-p^3t)} $$ with
\begin{eqnarray} deg(P^{(p)}_3(t)) &=& 2+2h^{(2,1)} \nonumber \\
deg(P^{(p)}_2(t)) &=& h^{(1,1)} \nonumber
\end{eqnarray}
This
follows from the fact that for non-toroidal Calabi-Yau threefolds
we have $b_1=0$.

({\bf d.}) Finally, Calabi-Yau fourfolds (with $\pi_1$ finite) are
of interest because of conjectured dualities involving Vafa's
F-theory \cite{v96}\cite{mv96}. The general cohomology structure
$${\small \matrix{
     &   &       &     &1       &     &      &    & \cr
     &   &       &0    &        &0    &      &    &   \cr
     &   &0      &     &h^{1,1}       &     &0      &    &   \cr
 &0  &       &h^{2,1}    &        &h^{2,1}    &       &0   &   \cr
 1 & &h^{3,1} &  &h^{2,2}    &     &h^{3,1}    &    &1,  \cr}} $$
gives rise to $$ Z(X_4/{\mathbb F}_p,t) =
\frac{P_3^{(p)}(t)}{(1-t)P_2^{(p)}(t)
P_4^{(p)}(t)P_6^{(p)}(t)(1-p^4t)} $$ with polynomials of the
degrees
\begin{eqnarray}
 {\rm deg} P_2^{(p)}(t) &=& {\rm deg} P_6^{(p)} = h^{(1,1)},
 \nonumber \\
 {\rm deg} P_3^{(p)}(t) &=& 2h^{2,1}, \nonumber \\
 {\rm deg} P_4^{(p)}(t) &=& 2+2h^{3,1}+h^{2,2}.
 \end{eqnarray}

{\bf 4.2} The relation between conformal field theories and
Calabi-Yau manifolds in the class of Brieskorn-Pham varieties, or
more general polynomials related to Kazama-Suzuki models
\cite{ls90c}, in general involves singular complete intersections
which must be resolved, and/or complete intersections in products
of weighted projective spaces. Such varieties have second Betti
numbers that can be quite large. The CFT-CY relation however must
hold already for the very simplest members in this class in terms
of cohomology. All the information necessary to establish such a
relation  must therefore be present already in these simple
examples.  Hence we can consider smooth Calabi-Yau hypersurfaces
in ordinary projective space $$\left\{\sum_i z_i^{n+1} =0 \right\}
~~\in ~~{\mathbb P}_{n+1}[n+2]$$ or in smooth hypersurfaces in
weighted projective space.
 Then
$$b_i(X_d) \in \{0,1\},~~i<d.$$

 For smooth threefolds embedded in weighted projective spaces this
 leads to a simplified zeta function of the form
 \begin{equation} Z(X/{\mathbb F}_p, t)=
\frac{P^{(p)}_3(t)}{(1-t)(1-pt)(1-p^2t)(1-p^3t)}
\label{zetacy3smooth}
\end{equation} with \begin{equation} deg(P^{(p)}_3(t)) =
2+2h^{(2,1)}\end{equation} and
\begin{equation} \sum_{i=1}^{b_3} \beta_i^{(3)} = 1+p+p^2+p^3
-\#(X/{\mathbb F}_p). \end{equation} Hence for smooth hypersurface
threefolds the only interesting information of the zeta function
is encoded in the polynomials $P_3^{(p)}(t)$.

\section{Hasse-Weil L-Function}

{\bf 5.1}
 We see from the rationality of the zeta function that the basic
 information
 of this quantity is parametrized by the cohomology of the variety.
 More precisely, one can show that  the $i$'th polynomial
 $P_i^{(p)}(t)$ is
 associated to the action induced by the Frobenius morphism on the
 $i$'th cohomology group ${\rm H}^i(X)$. In order to gain insight
 into the arithmetic information encoded in this Frobenius action
 it is useful to decompose the zeta function. This leads to the
 concept of a local L-function that is associated to the
 polynomials $P_i^{(p)}(t)$ via the following definition.

 Let $P_i^{(p)}(t)$ be the polynomials
 determined by the rational congruent zeta function over the field
 ${\mathbb F}_p$. The $i$'th L-function of the variety
 $X$ over ${\mathbb F}_p$ then is defined via
 \begin{equation} L^{(i)}(X/{\mathbb F}_p, s) =
 \frac{1}{P_i^{(p)}(p^{-s})}.\end{equation}

 Such L-functions are of interest for a number of reasons. One of
these is that often they can be modified by simple factors so that
after analytic continuation they (are conjectured to) satisfy some
type of functional equation.

{\bf 5.2} It was mentioned in the previous section that the
geometry/CFT relation must hold for the simplest type of
varieties, in particular those that do not need any kind of
resolution. When considering the cohomology of such simple
varieties in dimensions one through four it becomes clear that
independent of the dimension of the Calabi-Yau variety the
essential ingredient is provided by the cycles that span the
middle-dimensional (co)homology. Hence, even though in general the
cohomology can be fairly complex, in particular for Calabi-Yau
fourfolds,  the  only local factor in the L-function that is
relevant for the present discussion is $P_d^{(p)}(t)$ for $d={\rm
dim}_{{\mathbb C}}X$.

This suggests a natural generalization of the concept of the
Hasse-Weil L-function of an elliptic curve.
 Let $X$ be a Calabi-Yau $d$-fold with $h^{i,0}=0$ for $0<i<n-1$
 and denote by ${\mathfrak P}(X)$ the set of good prime numbers
 of $X$, i.e.
 those prime numbers for which the variety has good reduction over
 ${\mathbb F}_p$.
Then its associated Hasse-Weil L-function is defined as
\begin{equation} L_{\rm HW}(X,s) = \prod_{p\in {\mathfrak P}(X)}
\frac{1}{P^{(p)}_d\left(p^{-s}\right)} = \prod_{p\in {\mathfrak
P}(X)} \frac{1}{\prod_{j=1}^{b_3}
\left(1-\beta_j^{(d)}(p)p^{-s}\right)}.
\end{equation}  The existence of bad primes
complicates the whole theory considerably, but for our purposes we
can ignore the additional factors of the completed L-functions
that are induced by these bad primes.

As mentioned above, for smooth weighted CY hypersurfaces the
Hasse-Weil L-function then contains the complete arithmetic
information of the congruent zeta function. For reasons that will
become clear below it is of importance that the above Hasse-Weil
function can be related to a Hecke L-function, induced by Hecke
characters. The main virtue of such characters is that as the
simpler Dirichlet characters they are multiplicative maps. It is
this multiplicativity which is of essence for the present
framework.

\section{The Quintic}

 The Calabi-Yau variety which we will consider is the
quintic hypersurface in ordinary projective fourspace ${\mathbb
P}_4$. We denote the general $h^{2,1}=101$ complex dimensional
family of quintic hypersurfaces in projective fourspace ${\mathbb
P}_4$ by ${\mathbb P}_4[5]$ and consider the element defined by
\begin{equation} {\mathbb P}_4[5] \ni X= \left\{\sum_{i=0}^4
x_i^5=0\right\}.\end{equation} It follows from Lefshetz's
hyperplane theorem that the cohomology below the middle dimension
is inherited from the ambient space. Thus we have
$h^{1,0}=0=h^{0,1}$ and $h^{1,1}=1$, while $h^{2,1}=101$ follows
from counting monomials of degree five. Following Weil \cite{w49}
the zeta function is determined by (\ref{zetacy3smooth}), where
the numerator is given by the polynomial
$P_3^{(p)}(t)=\prod_{i=1}^{204} (1-\beta_i^{(3)}(p)t)$ which takes
the form  \begin{equation} P_3^{(p)}(t)= \prod_{\alpha \in
{\mathcal A}} \left(1-j_p(\alpha) t \right).\end{equation} This
expression involves the following ingredients. Define $\ell =
(5,p-1)$ and rational numbers $\alpha_i$ via $\ell \alpha_i\equiv
0({\rm mod}~1)$. The set ${\mathcal A}$ is defined as
\begin{equation} {\mathcal A}=
\{\alpha=(\alpha_0,...,\alpha_4) ~|~0<\alpha_i <1,~\ell
\alpha_i\equiv 0({\rm mod}~1), \sum_i \alpha_i=0({\rm
mod}~1)\}.\end{equation} Defining the characters
$\chi_{\alpha_i}\in \hat{{\mathbb F}_p}$ in the dual of ${\mathbb
F}_p$ as $\chi_{\alpha_i}(u_i)=exp(2\pi i\alpha_i s_i)$ with
$u_i=g^{s_i}$ for a generating element $g \in {\mathbb F}_p$, the
factor $j_p(\alpha)$ finally is determined as
\begin{equation} j_p(\alpha) =\frac{1}{p-1}\sum_{\sum_i u_i=0}
\prod_{i=0}^4 \chi_{\alpha_i}(u_i).\end{equation}

We thus see that the congruent zeta function leads to the
Hasse-Weil L-function associated to a Calabi-Yau threefold
\begin{equation} L_{\rm HW}(X,s) =
\prod_{p \in {\mathfrak P}(X)} \prod_{\alpha \in {\mathcal A}}
\left(1-\frac{j_p(\alpha)}{p^s}\right)^{-1},\end{equation}
ignoring the bad primes, which are irrelevant for our purposes.

\section{Algebraic Hecke Characters}

As mentioned in the introduction, the for the uninitiated
surprising aspect of the Hasse-Weil L-function is that it is
determined by another, a priori completely different, kind of
L-function. This second type of L-function is derived not from a
variety but from a number field. It is this possibility to
interpret the cohomological Hasse-Weil L-function as a field
theoretic L-function which establishes the connection that allows
to derive number fields $K$ from algebraic varieties $X$.

For the case at hand the type of L-function that is relevant is
that of a Hecke L-function determined by a Hecke character, more
precisely an algebraic Hecke character. Following Weil we will see
that the relevant field for our case is the extension ${\mathbb
Q}(\mu_m)$ of the rational integers ${\mathbb Q}$ by roots of
unity, generated by $\xi=e^{2\pi i/m}$ for some rational integer
$m$. It turns out that these fields fit in very nicely with the
conformal field theory point of view. In order to see how this
works this Section first describes the concept of Hecke characters
and then explains how the L-function fits into this framework.

There are many different definitions of algebraic Hecke
characters, depending on the precise number theoretic framework.
Originally this concept was introduced by Hecke \cite{h18} as
Gr\"ossencharaktere of an arbitrary algebraic number field. In the
following Deligne's adaptation of Weil's Gr\"ossencharaktere of
type $A_0$ is used \cite{d77}.

{\bf Definition.} {\it Let ${\mathcal O}_K \subset K$ be the ring
of integers of the number field $K$ and ${\mathfrak f} \subset
{\mathcal O}_K$ an ideal. Denote by ${\mathfrak I}(K)$ the set of
ideals of $K$. An algebraic Hecke character modulo ${\mathfrak f}$
is a multiplicative function $\chi$ defined on the ideals $
{\mathfrak I}(K)$ that are relatively prime to ${\mathfrak f}$ for
which the following condition holds. There is an element $\sum
n_{\sigma}\sigma \in {\mathbb Z}[{\rm Gal} (K/{\mathbb Q})]$ in
the integral group ring of the Galois group of the abelian
extension $K/{\mathbb Q}$ such that if $\alpha \in {\mathcal
O}_K$, $\alpha \equiv 1({\rm mod}~{\mathfrak f})$ then}
\begin{equation}
\chi((\alpha)) = \prod_{\sigma} \sigma(\alpha)^{n_{\sigma}}.
\end{equation}
{\it Furthermore, there is an integer $w>0$ such that
$n_{\sigma}+n_{{\bar \sigma}}=w$ for all $\sigma \in {\rm
Gal}(K/{\mathbb Q})$. This integer $w$ is called the weight of the
character $\chi$.}

Given any such character $\chi$ defined on the ideals of the
algebraic number field $K$ we can follow Hecke and consider a
generalization of the Dirichlet series via the L-function
\begin{equation} L(\chi,s)=\prod_{\stackrel{{\mathfrak p}
\subset {\mathcal O}_K}{{\mathfrak p} ~{\rm prime}}}
\frac{1}{1-\frac{\chi({\mathfrak p})}{{\rm N}{\mathfrak p}^s}} =
\sum_{{\mathfrak a} \subset {\mathcal O}_K} \frac{\chi({\mathfrak
a})}{{\rm N}{\mathfrak a}^s},\end{equation} where the sum runs
through all the ideals. Here ${\rm N}{\mathfrak p}$ denotes the
norm of the ideal ${\mathfrak p}$, which is defined as the number
of elements in ${\mathcal O}_K/{\mathfrak p}$. The norm is a
multiplicative function, hence can be extended to all ideals via
the prime ideal decomposition of a general ideal. If we can deduce
from the Hasse-Weil L-function the particular Hecke character(s)
involved we will be able to derive from the varieties
distinguished number field $K$.

Insight into the nature of number fields can be gained by
recognizing that for certain extensions $K$ of the rational number
${\mathbb Q}$ the higher Legendre symbols provide the characters
that enter the discussion above. Inspection then suggests that we
consider the power residue symbols of cyclotomic fields
$K={\mathbb Q}(\mu_m)$ with integer ring ${\mathcal O}_K =
{\mathbb Z}[\mu_m]$. The transition from the cyclotomic field to
finite fields is provided by the character which is determined for
any algebraic integer $x\in {\mathbb Z}[\mu_m]$ prime to $m$ by
the map
\begin{equation} \chi_{\bullet} (x): {\mathfrak I}_m({\mathcal O}_K)
\longrightarrow {\mathbb C}^*,\end{equation} which is defined on
ideals ${\mathfrak p}$ prime to $m$ by sending the prime ideal to
the $m$'th root of unity for which \begin{equation} {\mathfrak p}
~~\mapsto~~\chi_{\mathfrak p}(x)=x^{\frac{{\rm N}{\mathfrak
p}-1}{m}} ({\rm mod}~{\mathfrak p} ).\end{equation} Using these
characters one can define Jacobi-sums of rank $r$ for any fixed
element $a=(a_1,...,a_r)$ by setting
\begin{equation} J_a^{(r)}({\mathfrak p})=(-1)^{r+1} \sum_{\stackrel{u_i\in
{\mathcal O}_K/{\mathfrak p}}{\sum_i u_i=-1 ({\rm mod}~{\mathfrak
p})}} \chi_{\mathfrak p}(u_1)^{a_1}\cdots \chi_{\mathfrak
p}(u_r)^{a_r}
\end{equation} for prime ${\mathfrak p}$.  For non-prime ideals
${\mathfrak a} \subset {\mathcal O}_K$ the sum is generalized via
prime decomposition ${\mathfrak a} = \prod_i {\mathfrak p}_i$ and
multiplicativity $J_a({\mathfrak a})=\prod_i J_a({\mathfrak
p}_i)$. Hence we can interpret these Jacobi sums as a map
$J^{(r)}$ of rank $r$ \begin{equation} J^{(r)}: {\mathfrak
I}_m({\mathbb Z}[\mu_m]) \times ({\mathbb Z}/m{\mathbb Z})^r
\longrightarrow {\mathbb C}^*, \end{equation} where ${\mathfrak
I}_m$ denotes the ideals prime to $m$. For fixed ${\mathfrak p}$
such Jacobi sums define characters on the group $({\mathbb
Z}/m{\mathbb Z})^r$. It can be shown that for fixed $a\in
({\mathbb Z}/m{\mathbb Z})^r$ the Jacobi sum $J_a^{(r)}$ evaluated
at principal ideals $(x)$ for $x\equiv 1({\rm mod}~m^r)$ is of the
form $x^{S(a)}$, where
\begin{equation} S(a) = \sum_{\stackrel{(\ell,m)=1}{\ell~{\rm
mod}~m}}\left[\sum_{i=1}^r \left< \frac{\ell
a_i}{m}\right>\right]\sigma_{\ell}^{-1},\end{equation} where $<x>$
denotes the fractional part of $x$ and $[x]$ describes the integer
part of $x$.

We therefore see that the Hasse-Weil L-function is in fact a
product of functions each of which is determined by a Hecke
character defined by a Jacobi sum that is determined by a prime
ideal in the cyclotomic field ${\mathbb Q}(\mu_m)$. In the case of
the quintic hypersurface we derive in this way the fusion field
from the arithmetic structure of the defining variety reduced to a
finite field. To summarize, we have seen that the fusion field of
the underlying conformal field theory is precisely that number
field which is determined when the cohomological Hasse-Weil
L-function is interpreted as the Hecke L-function associated to an
algebraic number field.

\section{Quantum Dimensions}

{\bf 8.1} The discussion so far leads to the question whether
there is a natural field theoretic interpretation of the quantum
dimensions. The quantum dimensions are real, hence one needs to
consider the real field $${\mathbb Q}(\mu_m)^{+} \subset {\mathbb
Q}(\mu_m),$$ generated by $(\xi_m + \xi_m^{-1})$.

This field naturally emerges via the class number of the
cyclototomic field, which splits as $$ h({\mathbb Q}(\mu_m)) =
h^{+} ~h^{-} $$ where $$h^{+} = h({\mathbb Q}(\mu_m))^{+}.$$

The class number $h^{+}$ in turn singles out within the real field
${\mathbb Q}(\mu_m)^{+}$ the cyclotomic units $$ \theta_j =
\left|\frac{1-\xi_m^j}{1-\xi_m}\right| = \frac{\sin
\frac{j\pi}{m}}{\sin\frac{\pi}{m}}. $$ These numbers are precisely
the quantum dimensions!

The exploration of this has a long history, more than 150 years,
starting with Kummer.

{\bf Theorem.}[Kummer 1847, Sinnott 1978] {\it Let $m>2$, $m\neq
2~({\rm mod}~4)$ be the conductor of ${\mathbb Q} (\mu_m)$. Denote
by $U_c^{+}$ the subgroup spanned by the cyclotomic units within
the group $U^{+}$ of positive real units in ${\mathbb Q}(\mu_m)$.
Then $$h^{+} = 2^b [U^+:U_c^{+}],$$ where $b=0$ if the number $g$
of prime factors is unity, and $$b=2^{g-2} + 1 -g$$ if $g>1$.}

This provides an identification of the quantum dimensions.

{\bf 8.2} Furthermore, the class number $h^{+}$ of the maximal
real subfield ${\mathbb Q}(\mu_p)^{+}$ is constructed in part from
the real cyclotomic units. For $p$ an odd prime $$h^{+} =
\frac{2^{(p-3)/2}\Delta}{R},$$ where $R$ is the regulator of the
field, in general defined as the logarithmic image of a set of
fundamental units and $\Delta$ is a determinant constructed from
the quantum dimensions
 \begin{equation} \Delta ~=~\left|{\rm det}\left(
\sigma_j (\theta_k) \right)_{\stackrel{2\leq k\leq (p-1)/2}{0\leq
j\leq (p-3)/2}} \right|.\end{equation}

The regulator $R$ can be viewed as the volume of the logarithmic
image of a fundamental system of units. It was shown by Dirichlet
that the group $U$ of units in an algebraic number field of degree
$[K:{\mathbb Q}]=r_1+2r_2$ takes the form \begin{equation} U \cong
\mu \prod_{i=1}^{r_1+r_2-1} G_i,\end{equation} where $\mu$ is the
group of roots of unity, each $G_i$ is a group of infinite order,
and $r_1$ ($r_2$) denotes the number of real (complex) embeddings
of the field $K$. Hence every unit $u\in U$ can be written in the
form $u=\alpha \prod_{i=1}^r \epsilon_i$, where $\alpha \in \mu$
and $\{\epsilon_i\}_{i=1,...,r=r_1+r_2-1}$ is called a fundamental
system of units. It is useful to translate the multiplicative
structure of the units into an additive framework via the
regulator map
\begin{equation}
r: U \longrightarrow {\mathbb R}^{r_1+r_2}
\end{equation}
defined by
\begin{equation}
r(u) = \left(\ln|\rho_1(u)|,\dots,
\ln|\rho_{r_1}(u)|,\ln|\rho_{r_1+1}(u)|^2,\dots,
\ln|\rho_{r_1+r_2}(u)|^2\right),\end{equation} where
$\{\rho_i\}_{i=1,\dots,r_1}$ are the real embeddings and
$\{\rho_{r_1+j}\}_{j=1,\dots,r_2}$ are the complex embeddings of
$K$. The regulator $R$ then is defined as \begin{equation} R={\rm
det} \left( a_i \ln |\rho_i(\epsilon_j)|\right)_{\stackrel{1\leq
i\leq r_1+r_2}{1\leq j\leq r_1+r_2-1}}, \end{equation} with
$a_i=1$ for the real embeddings $i=1,...,r_1$, and $a_i=2$ for the
complex embeddings $i=r_1+1,...,r_1+r_2$. The regulator is
independent of the choice of the fundamental system of units.

{\bf 8.3} The results described above provide an example in which
the class number of an algebraic number field acquires physical
significance. This is not without precedence.
 Recently the class number of the fields of definition
of certain arithmetic black hole attractor varieties have been
interpreted as the number of U-duality classes of black holes with
the same area \cite{m98}. Here we see a further instance where the
class number of an algebraic number field acquires physical
meaning.

\section{Beilinson-Bloch Conjectures}

{\bf 9.1} The idea of analyzing the conformal field theory
structure of Calabi-Yau varieties in terms of their arithmetic
properties originally suggested itself \cite{s95} via some work by
Bloch and Schoen that was aimed at confirming some of the
conjectures of Beilinson \cite{ab85} and Bloch \cite{b84}. It was
shown in refs. \cite{b85} and \cite{s86} in the context of the
resolution of a nodal quintic threefold that the vanishing order
of the Hasse-Weil function at the central value determines the
rank of a twisted cohomology group, where the twist is determined
by a quadratic extension of the rational numbers.

More precisely, one has the following

{\bf Conjecture. (Swinnerton--Dyer, Beilinson, Bloch)} \hfill
\break
 Consider a complex variety $X$ and denote by {\rm CH}$^r(X)$
 the Chow group of codimension $r$ algebraic
cycles modulo rational equivalence. Let furthermore {\rm
CH}$^r(X)_0$ denote the subgroup of nullhomologous cycles and {\rm
L(H}$^{2r-1},s$) be the L$-$function associated to the cohomology
group {\rm H}$^{2r-1}$. Then $$ {\rm rk}~{\rm CH}^r(X)_0 = {\rm
ord}_{s=r}^{\rm zero} ~{\rm L(H}^{2r-1},s).$$

More precisely, these conjectures equate the vanishing order of
the L--function associated to an odd dimensional cohomology group
${\rm H}^{2r-1}(X_{{\bar K}},{\mathbb Q}_{\ell})$ of a smooth
projective variety $X$ over a number field $K$ at the point $s=r$
to the rank of the group $$A^r(X)=Ker\left({\rm CH}^r(X_K)\otimes
{\mathbb Q} \longrightarrow H^{2r}(X_{{\bar K}},{\mathbb
Q}_{\ell}(r))\right),$$  where CH$^r(X)$ is the Chow group of
codimension $r$ cycles on $X$ defined over $K$ modulo rational
equivalence. Then $$ {\rm rk~A}(X) = {\rm ord}_{s=r}~{\rm
L(H}^{2r-1},s). $$
 \vskip .2truein

{\bf 9.2 Example.} A tantalizing check for the Beilinson-Bloch
conjectures was considered by Bloch \cite{b85} and Schoen
\cite{s86}. Their starting point is the deformed quintic at the
conifold point
\begin{equation} Y_0 = \{z_0^5 + \cdots + z_4^5 - 5\prod_{i=0}^4
z_i =0\} ~\subset {\mathbb P}_4.\end{equation} This variety fails
to transverse at 125 nodal points $$
\left\{(\xi_0,...,\xi_4)~|~\xi_i^5=1,~\prod_{i=1}^4 =1\right\}. $$
The resolution was described in detail by Schoen \cite{s86} and
was found to lead to a rigid manifold $Y$ with Euler number
$\chi(Y)=50$ and Hodge diamond

\begin{footnotesize}
$$ Y:~~~~
\begin{tabular}{c c c c c c c}
  &  &          &1         &              &   &   \\
  &  &0         &          &0             &   &   \\
  &0 &          &25        &              &0  &   \\
1 &  &0         &          &0             &   &1  \\
  &0 &          &25        &              &0  &   \\
  &  &0         &          &0             &   &   \\
  &  &          &1         &              &   &   \\
\end{tabular}
$$
\end{footnotesize}

The interesting thing now is that even though this variety is
rigid over ${\mathbb C}$, the intermediate cohomology does not
vanish when one considers cohomology with weird fields. Namely, it
was shown by Bloch and Schoen that the following holds.

{\bf Theorem.}[Bloch] \hfill \break 1. $L(H^3(Y),2)\neq 0$. \hfill
\break 2. {\it If $\rho$ is the quadratic character associated to
${\mathbb Q}(\sqrt{5})/{\mathbb Q}$
   then $L(\rho \otimes H^3(Y),s)$ vanishes to order 1 at $s=2$.}

It follows from standard conjectures and the theorem that the
Griffiths group in codimension 2 is of rank 0 over ${\mathbb Q}$
but has rank 1 over ${\mathbb Q}(\sqrt{5})$.  Assuming the purity
conjecture for \'etale cohomology the cycle implied by this result
has infinite order in the Griffiths group. In this way a nonzero
map is obtained $$ A^2(X_{{\mathbb Q}(\sqrt{5})}) \longrightarrow
{\mathbb Z}/3{\mathbb Z}.$$

The important point of all this is that these arithmetic
considerations single out the quadratic field ${\mathbb
Q}(\sqrt{5})$, which is precisely the field generated by quantum
dimensions of the Gepner model.

This is the good news. The bad news is that the significance of
this result is not obvious.

 {\bf (1)} The starting point is a deformed quintic, not the
Fermat quintic. Only the latter is known to be exactly solvable.
Deformed conformal field theories are still poorly understood,
hence one could question the significance of this result. One
argument in favor of its relevance is that the quantum dimensions
are indicative of cohomology, hence deformations per se should not
produce dramatic behavior (unless one hits a singular point).

{\bf (2)} Worse, the starting point is nodal, precisely when
drastic things do happen, and the computation is done for the
resolved variety.

In light of these difficulties the obvious question arises whether
there are other varieties which feature a similar behavior and
which can be related to exactly solvable field theories. In this
context a conjecture by Yui is rather interesting.

{\bf Conjecture.}[Yui]{\it For any Brieskorn-Pham threefold
defined over ${\mathbb Q}$ there exists a character $\rho$
associated to the cyclotomic field such that} $${\rm
ord}_{s=2}~{\rm L}(\rho\otimes {\rm H}^3(X_{{\mathbb Q}}),2) \geq
1.$$


\end{document}